# Mathematical Models and Formal Approaches for Caching Strategies Optimisation

## Vitaly.O. Groppen

*(North-Caucasian Institute of Mining and Metallurgy, Vladikavkaz, Russia)*

An optimization of caching strategies is proposed as a formal approach allowing us a more efficient use of two-level computer memory. This approach is based on a set of mathematical models and a set of theorems, permitting analytical determination of optimal caching strategy operating characteristics, thus minimising time loss caused by optimal strategy search.

## 1. Introduction

This work is intended to help a user to use more efficiently computer multi-level memory [1, 6, 8]. It is based on optimal data allocation and caching strategy resulting, as a rule, in minimisation of external memory calls number [6-8]. Selection of caching strategy is interpreted as a selection of the size of each cache buffer in the main memory used for external data file caching. Proposed formal approach permits an analytical determination of optimal caching strategy for any composition, consisting of:
- computer equipped with two-level memory,
- problem being solved, with fixed size of each data file.

Mathematical models used ignore integer nature of variables reflecting size of each cache buffer. The following symbols and assumptions are used:

## 2. Symbols and assumptions used

$H$ - external data files number;
$N_i$ - upper bound of number of calls of i-th file ($i = 1, 2,..., H$);
$W_i$ - size of i-th external file;
$U_i$ - size of cache buffer in the main memory, being used for i-th user's file processing;
$\vec{U}$ - vector of variables with i-th component equal to $U_i$, ($i = 1, 2,..., H$);
$V$ - size of free main memory.

We further assume that:
- ratio $W_i/U_i$ ($i=1, 2, \ldots H$) is always integer;
- $U_i$ variables are continuous for all the values of index "i".

## 3. Problems statement

Minimising external memory calls number while processing any external data file, we, as a rule, use caching of each such a file by a cache buffer, the latter being allocated either in the special cache-memory, or in the main memory of the computer used [1, 8]. Problem statement for the last caching strategy in a general case can be presented as a multi-criteria problem [2, 3] - for each external file external memory calls number is minimised under the following conditions:
1. For each i-th external file caching is used i-th cache-buffer.

2. Memory size of each cache-buffer serves as a tool of optimisation.
3. Memory size of each external user's file:
   - does not depend on time;
   - is known *a priory*.
4. Total size of all cache-buffers does not exceed free main memory size V of the computer used.

Keeping in mind all these conditions we can suggest the following formal problem statement:

(1) $$\begin{cases} \forall i, \ f_i = W_i/U_i \to \min; \\ \sum_i U_i = V; \\ \forall i, \ 1 \leq U_i \leq W_i. \end{cases}$$

As a rule the search of multi-criteria problem solution is based on Pareto optima principle [3]: vector of variables is optimal if improvement of values of any subset of goal functions is always accompanied by a change for the worse of another goal functions subset. It is obvious, that system (1) Pareto-optimal vector of variables is simultaneously optimal for the problem with maximised inverse goal functions:

(2) $$\begin{cases} \forall i, \ \varphi_i = U_i/W_i \to \max; \\ \sum_i U_i = V; \\ \forall i, \ 1 \leq U_i \leq W_i. \end{cases}$$

Moreover, as $W_i$ values are constants (i = 1, 2,.., H), system (2) Pareto-optimal vector of variables is also optimal for the following system:

(3) $$\begin{cases} \forall i : \psi_i = U_i \to \max; \\ \sum_i U_i = V; \\ \forall i, \ 1 \leq U_i \leq W_i. \end{cases}$$

An imperfection of Pareto-optima definition seems to be its' low selectivity: power of set of optimal solutions is often commensurable with a set of all feasible plans [4]. To avoid this, a multi-criteria problem is usually transformed into a problem with a single goal function by means of certain additional information, not available in the original problem statement. Example of such approach is creation of super-criterion generated as a sum of weighted goal functions of original problem [2]. Using number of calls of each i-th file $N_i$ (i=1,2,…, H) as a natural weight of system (1) i-th goal function, we can present the above mentioned additional information known *a priory,* and hence minimising the total number of external memory calls number, system (1) is substituted by the following one [7]:

$$(4) \quad \begin{cases} F_1 = \sum_i \dfrac{W_i}{U_i} N_i y_i \to \min; \\ \forall i, \ y_i = signum(W_i - U_i); \\ \sum_i U_i = V; \\ \forall i, \ 1 \leq U_i \leq W_i. \end{cases}$$

If we aim at minimization of upper bound of external memory calls number for any file, either system (1) goal functions are replaced by the following objective:

$$(5) \quad F_2 = \max_i \dfrac{W_i}{U_i} \to \min,$$

or system (4) goal function is replaced by (6):

$$(6) \quad F_3 = \max_i \left\{ N_i y_i \dfrac{W_i}{U_i} \right\} \to \min.$$

Any intelligent database (IDB), using (4) and (6) for external memory calls number minimization via optimal values of cache-buffers $U_i$ determination, must permanently collect and periodically renew values $N_i$, thus losing more time for IDB self-service [6]. Elimination of these procedures is based on the presented below approach. Transformation principles for a multi-criteria problem into a problem with a single goal function that does not need any additional condition with the use of comparison standards are proposed in [4]. In the studied above case these standards are connected with combinations of the best (comparison standard "a") and the worst (comparison standard "b") values of system (1) goal functions regardless of the existence of corresponding feasible vectors of variables $\vec{U}$.

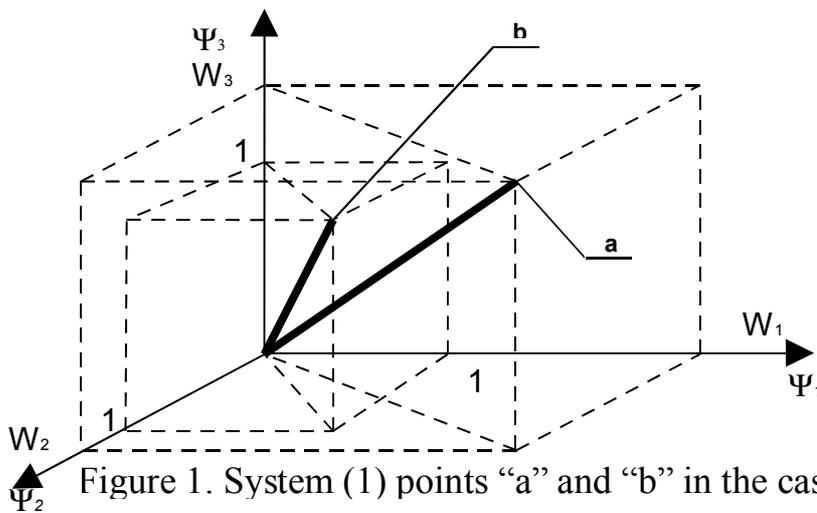

Figure 1. System (1) points "a" and "b" in the case of

It is easy to see that in H-dimensional criteria space problem (3) comparison standard "b" is reflected by a point, corresponding to the unit criteria vector $\vec{K} = \{1,1,......,1\}$, whereas standard "a" is presented by a point corresponding to vector $\vec{W} = \{W_1, W_2, ..., W_n\}$ (Figure 1).

Using comparison standard "a", and proved in [4] theorems, searching of system (3) solution can be replaced by solving the following problem:

(7)
$$\begin{cases} F_4 = \sum_i [W_i - U_i]^2 \to \min; \\ \sum_i U_i = V; \\ \forall i,\ 1 \le U_i \le W_i. \end{cases}$$

In other words (7) reflects searching of permissible vector $\vec{U}$, corresponding to some point "c", outside the unit hypercube, simultaneously being the nearest to the point "a" in H-dimensional criteria-space (Figure 1).

The features of problems (1) – (7) optimal solutions are discussed below.

## 4. Features of optimal solutions

Several theorems below reflect useful features of problems (2) – (7) optimal solutions, permitting analytical determination of optimal caching strategy, thus minimising a time loss caused by cache-buffers optimization.

**Theorem 1.** The lower bound of system (4) goal function value is presented by the following vector of variables [1, 6, 7]:

(8) $\quad \forall i:\ U_i = V \dfrac{\sqrt{N_i W_i}}{\sum\limits_{j=1}^{H} \sqrt{N_j W_j}}.$

For problem resulting in substitution of system (4) goal function by goal function (5), the following is true [8]:

**Theorem 2.** Goal function $F_2$ value lower bound is presented by the following vector of variables:

(9) $\quad \forall i:\ U_i = V \dfrac{W_i}{\sum\limits_{j=1}^{H} W_j}.$

System (9) results in $F_2$ minimal value equal to $F_2'$:

(10) $\quad F_2' = \dfrac{1}{V} \sum\limits_{j=1}^{H} W_j.$

Similar result for the case, when system (4) goal function $F_1$ is replaced by $F_3$, is reflected in [8]:

**Theorem 3.** For formal statement of any problem, combining goal function (6) and system (4) restrictions, optimal vector of variables and corresponding goal function $F_3$ value are equal to:

(11)
$$\begin{cases} \forall i,\ U_i = V \dfrac{N_i W_i}{\sum\limits_{j=1}^{H} N_j W_j}; \\ F_3 = \dfrac{1}{V} \sum\limits_{j=1}^{H} N_j W_j. \end{cases}$$

It is easy to see that, if $\forall i,\ N_i = 1$, then systems (9) and (11) coincide.

System (7) optimal solution, corresponding to the nearest to standard "a" point "c" (Figure 1) with permissible vector of variables $\vec{U}$, can be expressed by the following theorem:

**Theorem 4.** System (7) optimal vector of variables is equal to:

$$(12) \quad \forall i, \ U_i = \begin{cases} W_i, & \text{if } V \geq \sum_j W_j; \\ 1, & \text{if } H = V; \\ W_i - \dfrac{\sum_j W_j - V}{H}, & \text{if } V < \sum_j W_j. \end{cases}$$

**Proof of theorem 4.**

Transforming (7) into Lagrange function L and equating $\dfrac{\partial L}{\partial U_j}$, $j = 1,2,\ldots,H$, $\dfrac{\partial L}{\partial \lambda}$ to zero, we obtain the following system:

$$(13) \quad \begin{cases} \forall j: -2(W_j - U_j) + \lambda = 0; \\ V - \sum_j U_j = 0. \end{cases}$$

System (13) solution looks like:

$$(14) \quad \begin{cases} \lambda = -2 \dfrac{\sum_j W_j - V}{H}; \\ \forall i, \ U_i = \begin{cases} W_i, & \text{if } V \geq \sum_j W_j; \\ W_i - \dfrac{\sum_j W_j - V}{H}, & \text{if } V < \sum_j W_j. \end{cases} \end{cases}$$

Eliminating in (14) the first equality and keeping in mind that $\forall i, \ 1 \leq U_i \leq W_i$, we transform this system into (12), what was to be proved.

### Conclusion

Though general model (1) was developed and used for the search of optimal caching strategy, the latter often can not be determined uniquely due to the numerous different approaches used for problem (1) solving. It results in a set of different optimal caching strategies caused by different Pareto-optimal $\bar{U}$ vectors. Selection of such an approach depends on user's goals and mentality, but in any case its positive aspect due to theorems 1 – 4 above is in the possibility of optimal caching strategy analytical determination, eliminating usage of iterative procedures. It permits to ignore time loss caused by optimal strategy determination, thus improving efficiency of caching optimization. Experimental analyses of proposed above approach confirms its efficiency for a number of data domains [9].

### References

1. Vitaly O. Groppen, Vladimir V. Mazin. Efficient strategy for computer main and external memory interconnection. Proceedings of the XI All-Union Workshop on Control Problems. Tashkent, 1989, p. 202 (Russ).
2. Oleg I. Larichev. Solution Making: Theory and Methods. Logos, 2002 (Russ).
3. Pareto V. Cours d'Economie Politique. Lausanne: Houge, 1889.

4. Vitaly O. Groppen. Quality rating principles for multi-criteria problems solutions based on comparison standards. Proceedings of the fifth International Workshop "Steady development of mountain territories: problems and perspectives of science and education integration", Vladikavkaz, 2004, pp. 572 – 580 (Russ).
5. Vitaly O. Groppen. New Solution Principles of Multi-Criteria Problems Based on Comparison Standards, 2005, http://www.arxiv.org/ftp/math/papers/0501/050157.pdf
6. Vitaly O. Groppen. Expert systems for computer memory optimal control. Proceedings of the IFAC Workshop on Safety, Reliability and Applications of Emerging Intelligent Control Technologies. Hong Kong, 12-14 December 1994, pp. 226-228.
7. Vitaly O. Groppen. Smart Computing. Publishing House "Terek", Vladikavkaz, 2004, 104 p.
8. Igor Kopilov. Optimal Strategy of Disk Arrays Usage and Relational Database Files Decomposition. Proceedings of the Conference on New Information Technologies in Science, Education, Economy, 10-11 June, 2002, Vladikavkaz, Russia, pp.76-77 (Russ).
9. Igor V. Kopilov. Smart Computing: Optimal Data Exchange Between the Main and External Memory. Proceedings of the fifth Conference on Evolutionary Methods on Design, Optimization and Control with Applications to Industrial and Social Problems. Barcelona, September 15 – 17, 2003, p. 135.